\magnification=1200
\def\today{\ifcase\month\or
  January\or February\or March\or April\or May\or June\or
  July\or August\or September\or October\or November\or December\fi
  \space\number\day,\space\number\year}
\def\Z{{\mathchoice{{\bf Z}}{{\bf Z}}{{\rm Z}}{{\rm Z}}}}
\def\naturals{{\bf N}}
\def\reals{{\bf R}}
\def\given{\, \big|\,}

\def\half{{{1}\over{2}}}
\def\pr{{\rm I\!\!P}}
\def\bl{\bigl [}\def\br{\bigr ]}\def\cP{{\cal P}}
\def\cT{{\cal T}}\def\cS{{\cal S}}\def\cF{{\cal F}}
\def\cR{{\cal R}}\def\cE{{\cal E}}
\def\bm{B^{\rm min}}\def\o{\omega}\def\O{\Omega}
\def\cI{{\cal I}}\def\cC{{\cal C}}
\def\ZZ{\Z\cup\Z+\half}
%
%
%
\def\introd{1}
%
%
\def\introdassone{1}
\def\introdasstwo{2}
\def\introdremarkone{1}
\def\introdremarktwo{2}
%
%
\def\dynamicsxi{1}
\def\dynamicseta{2}
\def\relation{3}
\def\otherrelation{4}

%
%
\def\invar{2}
%
%
\def\invartheoremone{1}
\def\invartheoremtwo{2}
\def\invartheoremthree{3}
\def\invartheoremfour{4}
\def\invarassone{1}
\def\invarremarkone{1}
\def\invarremarktwo{2}
\def\invarremarkfour{3}
\def\decomposition{1}
\def\bum{2}
\def\propone{3}
\def\proptwo{4}
\def\primeone{5}
\def\primetwo{6}
%
%
\def\hydrod{3}
%
%
%
\def\hydrodtheoremone{1}  
\def\hydrodtheoremtwo{2}
\def\hydrodthethree{3}
\def\hydrodtheoremthree{3}

\def\hydrodexampleone{1}
\def\hydrodexampletwo{2}
%
%
\def\efen{1}
\def\state{2}
\def\frombilin{3}
\def\closeness{4}

\def\bernoulli{6}
\def\bernoulliautomat{7}
\def\inapplication{8}
\def\true{9}
%
\def\rate{4}
%
%
\def\ratetheorem{1}
\def\rateassertiontwo{1}
\def\rateassertionthree{2}
\def\rateremarkzero{1}
\def\rateremarkone{2}
%
%

\def\intheoremrate{2}
\def\tutu{3}
\def\titi{4}
\def\initial{5}
\def\notrandomwalk{6}
%
%

%

%
%

%
%

%
%
\def\bennaim{BRL}
\def\wolfone{W1}
\def\wolftwo{W2}
\def\bramson{BN}

\def\elskens{EF}
\def\cyclic{Fi}

\def\krugspohn{KS}
\def\belfer{BF}
\def\ermakov{E}
\def\ferravi{FR}
\def\feller{F}
\def\krapivsky{KSL}
\def\bilin{Bi}
\def\moving{B}
\def\karatzas{KaS}
\def\japonez{Y}
\def\review{Nl}
\def\computersim{N}
\def\liggett{L}
\def\informatic{CST}
\def\bertocci{Ber}
\def\krugspohntwo{KS2}
%
%

\centerline{\bf Invariant Measures and Convergence for Cellular Automaton 184}
\centerline{\bf and Related Processes.}
\vskip2mm
\centerline{Vladimir Belitsky\dag \ \ and 
Pablo A. Ferrari\footnote\dag{Instituto de
  Matem\'atica e Estat\'\i stica, Universidade de S\~ao Paulo, Cidade
  Universitaria, Rua de Mata\~ao, 1010, 05508-900, SP,
  Brazil. pablo@ime.usp.br, vbel1@ime.usp.br}}
\vskip1mm
\centerline{\it Universidade de S\~ao Paulo}
\vskip3mm
%
%
%
\noindent{\bf Abstract.}\ \ For a class of one-dimensional cellular automata, 
we review and complete the
characterization of the invariant measures (in particular, all invariant
phase separation measures), the rate of convergence to
equilibrium,  and the derivation of 
the hydrodynamic limit. The most widely known representatives of
this class of automata are: Automaton 184 from the classification of 
S. Wolfram [\wolfone], an annihilating particle system and a 
surface growth model.

\vskip4mm
\settabs \+{\sl Key words and phrases:}\ \ &\cr
\+{\sl Key words and phrases:}\ \ & Cellular automata, automaton 184,
ballistic annihilation,\cr
\+ &  annihilating deterministic motions, surface growth,\cr
\+ &invariant measures, hydrodynamic limits, rate of convergence\cr
\+ & to equilibrium.\cr 
\vskip2mm

\noindent{\sl AMS Classification Numbers:}\ \ 60K35
\vskip2mm


\noindent{\bf 1. INTRODUCTION}
\vskip2mm

{\bf Cellular Automaton 184} (CA 184) is a discrete time process
with state space $\{0,1\}^{\Z}$ and the following evolution rule:
if $\eta\in \{0,1\}^{\Z}$ is the state at time $n$ then the state 
$\eta^\prime$  at time $n+1$ is defined by
$$
\eta^\prime(x):=\cases {1,  & if $\eta(x)=\eta(x+1)=1$\cr
                    1,  & if $\eta(x)=1-\eta(x-1)=0$\cr
                    0   & otherwise\cr}\kern3em \forall\,\, x\in \Z
\eqno(\introd.\dynamicsxi)
$$
Here $\eta(x)$ denotes the value of $\eta:\Z\rightarrow \{0,1\}$ at the 
coordinate $x$. The dynamics of CA~184 will be denoted by the operator 
$C:\{0,1\}^{\Z}\rightarrow\{0,1\}^{\Z}$ defined by $C\eta=\eta^\prime$. Thus,
formally, CA 184 refers to a sequence $\{\eta_n, n\in \naturals\}$ 
such that $\eta_{n+1}=C\eta_n\, ,\, \forall n$. 

CA 184 models deterministic motions of identical particles on $\Z$ that obey 
the following rules: there may be at most one particle per site, and at
each integer time each particle inquires whether the site of $\Z$ to the right
of its current position is empty of another particle, and if it is so then 
it instantaneously jumps to this site. These rules follow immediately from 
(\introd.\dynamicsxi), if one interprets $\eta(x)=1/0$ as `` presence/absence 
of a particle at the site $x\in\Z$ in $\eta$''.

The number ``184'' in the name of this process is due to  the
classification of  Wolfram [\wolfone] (see also [\wolftwo]) of a class of
cellular automata. CA 184 has been 
used to model traffic ([\computersim], [\review] and references therein) and
to classifies densities in binary strings ([\informatic]). CA 184 has
two ``stochastic'' counterparts. One is the Totally Asymmetric Simple 
Exclusion  Process (TASEP). The second one has the same evolution rule as 
CA~184 but with ``a 
noise'' that is introduced by setting that each particle that can jump will do 
so with probability $p$ independently of anything else; we thus call
this process ``CA 184 with noise'' (CA\&N). The invariant 
measures for TASEP and CA\&N have been characterized in [\liggett] 
and [\japonez], respectively (curiously, only the case $p\leq 1/2$ was 
studied in [\japonez], but we believe that a similar technique allows to 
extend the results to any $p$) but the approach employed there 
(which is a stochastic coupling) does not apply to CA 184 (because of the lack
of stochasticity in its dynamics). 
\vskip2mm

{\bf Annihilating particle system} studied here can be also found in
the literature under the names Ballistic Annihilation (BA) and Annihilating 
Deterministic Motions. We shall adopt here the name BA.  BA is a discrete time
process with state space $\{-1,0,1\}^{\Z}$ and the following dynamics: if 
$\zeta$ is the state at time $n$ then the state $\zeta^\prime$ at time $n+1$ 
satisfies
$$
\zeta^\prime(x)=\cases{ 1, & if $\zeta(x-1)=1$ and neither $\zeta(x)=-1$ nor
both $\zeta(x)=0, \zeta(x+1)=-1$\cr   
-1, & if $\zeta(x+1)=-1$ and neither $\zeta(x)= 1$ nor both $\zeta(x)=0,
\zeta(x-1)=1$\cr 
                    0, & otherwise\cr}
\eqno(\introd.\dynamicseta)
$$
for all $x\in \Z$. Here $\zeta(x)$ denotes the value of $\zeta:\Z\rightarrow 
\{-1,0,1\}$ at the coordinate $x$. The dynamics of BA will be denoted by the 
operator $A:\{-1,0,1\}^{\Z}\rightarrow\{-1,0,1\}^{\Z}$ defined by 
$A\zeta=\zeta^\prime$. Thus, formally, BA refers to a sequence $\{\zeta_n,
n\in  \naturals\}$ such that 
$\zeta_{n+1}=A\zeta_n\, ,\, \forall n$.

BA may be also interpreted in terms of particles. We shall call them 
A-particles in order to distinguish  them from those that move in CA 184. 
The values $0,1, -1$ of  $\zeta(x)$ are interpreted by  saying that the 
site $x$ is respectively, free of an A-particle, contains an A-particle with 
velocity  1, and contains an A-particle with 
velocity $-1$. In the terms of particles, the dynamics 
of BA acquires the following interpretation: each A-particle  moves
along $\reals$ with its velocity (going in the direction to $-\infty$ 
($+\infty$), if the velocity is negative (positive, resp.)) and annihilates 
when meets another A-particle; upon annihilation both A-particles disappear
from the system forever. To link this interpretation with 
(\introd.\dynamicseta), it is necessary to note that if A-particles occupied
only the sites of $\Z$ then after a unit of time has elapsed the survived
A-particles are again solely at the sites of $\Z$.

BA is a natural model for the reaction $A+B\rightarrow inert$. For an 
extensive list of references in this respect, we refer a reader to
[\elskens]. A modification of BA in which the set of possible particle
velocities is larger than $\{-1, 1\}$ has been considered 
in the works [\bennaim], [\krapivsky]. We also mention here the work 
[\ermakov] which studies a coalescing particle system employing 
a very simple relation between it and  BA. This coalescing particle system is
a process in which particles move as in BA but coalesce instead of annihilate, 
upon collision;
coalesced particles choose then a new velocity from $\{-1, 1\}$ with equal
probabilities.
\vskip2mm

{\bf Surface Growth Model} (SG) is a discrete time  process with 
state space $\cR$, the space of piecewise linear functions (from $\reals$ 
to $\reals$) that have  slope either $+1$, $-1$, or $0$ between any two
consequent integer abscissas, and the following dynamics:
If $f(\cdot)$ is the state at time $n$ then $\hat f(\cdot)\in \cR$, the
state at time $n+1$, is obtained from $f(\cdot)$ by the following rule:
if $x\in \Z$ is a local minimum of $f$, then the values of $\hat f$ at
$[x-1, x+1]$ are obtained by reflecting the graph of $f$ around the
straight line which links the points $(x-1, f(x-1))$ and $(x+1, f(x+1))$; at
any point $y\in \reals$ which does not have a local minimum of $f$ in its
neighborhood of radius $1$, we set $\hat f(y)=f(y)$. SG will be the name of a
sequence $\{f_n(\cdot), n\in \naturals\}$ such that $f_{n+1}(\cdot)=\hat
f_{n}(\cdot)\, ,\, \forall n$.

Imagine $f(\cdot)$ as the
surface of a two-dimensional solid above some reference horizontal line.
Imagine then that diamond shaped particles of the side length 1 are thrown
on this solid, and those of them that fall in local minima of $f(\cdot)$ stick
to the solid while others disappear. The surface of the new solid will be then
what we have defined as $\hat f(\cdot)$. This justifies the name of the
process just introduced. We observe that it is also known under the name 
polynuclear growth model (PNG) and (one of its modifications) was first
studied almost 30 years ago  in [\bertocci] (see also [\krugspohntwo] for 
a list  of more recent works).
\vskip2mm

{\bf Equivalence of dynamics} of two processes means here an  
existence of a transformation $T$ from the state space of one process to 
that of another one such that if $T\eta=\zeta$ then 
$T\eta^\prime=\zeta^\prime$, where $\eta^\prime$ and $\zeta^\prime$ are the 
states of these processes at time $n+1$
given their states at time $n$ were $\eta$ and $\zeta$ respectively.
The fact that the dynamics of CA 184, BA and SG are equivalent has been known
since at least the work [\krugspohn].  We present it in 
Assertions~\introd.\introdassone\ and~\introd.\introdasstwo\ below. The proofs
are straightforward and thus, omitted.
\vskip2mm

\noindent{\bf Assertion~\introd.\introdassone.}\ (equivalence between CA~184
and BA)\ \
{\sl Define $T_{184, BA}:\{0,1\}^{\Z}\rightarrow\{-1,0,1\}^{\Z}$ by
$$
(T_{184, BA}\eta)(i)= 1-\eta (i)-\eta(i+1), \kern2em i\in \Z
\eqno(\introd.\relation)
$$                 
Then, if for some $\eta$ and $\zeta$, $T_{184, BA}(\eta)=\zeta$ then
$T_{184, BA}(\eta^\prime)=\zeta^\prime$, where $\eta^\prime$ relates to $\eta$ 
via (\introd.\dynamicsxi), and $\zeta^\prime$ relates to $\zeta$ via
(\introd.\dynamicseta).}
\vskip2mm

\noindent{\bf Remark~\introd.\introdremarkone.}\ $T_{184, BA}$ is
an injection but not a surjection. Indeed, if any two 
consequent sites of $\Z$ are occupied by A-particles
with opposite  velocities in a configuration $\zeta\in \{-1,0,1\}^{\Z}$ then
there does not exist a $\eta\in \{0,1\}^{\Z}$ such that $T_{184,
BA}\eta=\zeta$. As a consequence, a certain care is necessary, when CA 184 is
studied with the help of BA (the reader will see it in the proofs of
Theorems~\invar.\invartheoremtwo\ and \invar.\invartheoremthree).
\vskip2mm

\noindent{\bf Assertion~\introd.\introdasstwo.}\ (equivalence between BA and
SG)\ \ {\sl Define  $T_{BA, SG}:\{-1,0,1\}^{\Z}\rightarrow \cR$ of 
the following form
$$
\eqalign{(T_{BA, SG}\zeta)(n)-(T_{BA, SG}\zeta)(n-1)&=\zeta(n), \,\, 
n \in \Z\cr
(T_{BA, SG}\zeta)(0) & =0 \cr} \eqno(\introd.\otherrelation)
$$
(certainly, it is sufficient to determine $f\in \cR$ on $\Z$).
Then, if  $T_{BA,SG}\zeta=f+\alpha$ for some $\zeta\in \{-1,0,1\}^{\Z}, 
f\in\cR$ and $\alpha\in \reals$, then
$T_{BA, SG}\zeta^\prime=\hat f+\beta$, where $\zeta^\prime$ relates to 
$\zeta$ via (\introd.\dynamicseta), $\hat f$ relates to $f$ via the dynamics
of SG and $\beta$ is a real number whose value depends on $\alpha$ and
$f$.(Above, the relation $g=h+\alpha$ for functions $g$ and $h$ and a number
$\alpha$ means $g(x)=h(x)+\alpha\, \forall x\in \reals$.)}
\vskip2mm

\noindent{\bf Remark~\introd.\introdremarktwo.}\  The actual values of the
constants $\alpha$ and $\beta$ are ignored because they will be of no 
importance, when SG is employed 
to investigate BA and CA 184, and also when BA or CA 184 is employed to
investigate phenomena related to modification of the shape of the surface in
SG. In this respect we observe that in SG, the surface both grows and changes
its shape, but only the phenomena related to the shape change will be of
interest to us.
\vskip2mm

{\bf The contents of this paper:}\ \ 
Section~\invar\ characterizes the invariant measures for the dynamics
of the processes defined above.
This is done with a help of a quite elementary property of BA presented in
 Assertion~\invar.\invarassone; it states at any time $n$ the
distance  between two diverging subsequent particles 
(i.e., no any other particle in between) cannot be less than
$2n+1$.  (We note that this property has been noticed in [\informatic] where 
it was employed to show that CA 184 on a ring of size $N$ will come to its 
invariant state by time not greater than  $N/2$.) The set of the invariant 
measures which are also translation invariant is then easily classified 
(Remark~\invar.\invarremarktwo). For BA, an extremal measure of this set 
is supported either by  configurations that have only particles with 
velocity $+1$ or by those that have  only particles with velocity $-1$, 
moreover, under such a measure the positions of these particles
form a stationary ergodic process on $\Z$. Applying then
Assertions~\introd.\introdassone\ and \introd.\introdasstwo\ to this result 
one easily characterizes the translation invariant measures that are
invariant for CA~184 and SG. For  extremal translation invariant measures 
of CA~184 we calculate the flux of particles through zero. The curve of the 
flux (as a function of the density of particles $\rho\in [0,1]$) is 
$1/2-|1/2-\rho|$; we observe that it is known to be $\rho(1-\rho)$ for 
TASEP. 
 
The invariant measures of BA which are not translation invariant are called 
``phase separation'' measures. These measures are supported by
configurations in which all particles with velocity $+1$ are situated to the 
left of all those whose velocity is $-1$. Our characterization  of the 
phase separation measures is obtained
via introducing in BA an extra ``second class'' particle in the way such that
its paths record the positions of particles in BA, and studying the properties
of these paths. The second class particle is always between the rightmost
particle with velocity $+1$ and the leftmost particle with velocity $-1$. It
moves together with the latter till it is annihilated by the former.
At this moment, the second class particle changes its velocity and 
starts moving to the right until it meets another particle with velocity
$-1$.  The definition of the second class particle yields the following
relations of its path to particle positions in BA:
each time interval the second class particle moves to the right (respectively,
left), is equal to the distance between two subsequent particles with velocity
$-1$ (respectively, $+1$). As for the relative positions of particles with
velocity $-1$ to those with velocity $+1$, it must be such that the their
annihilation times form a point process identical to that formed by the times
of change of velocity of the second class particle. These relations are
employed in Theorem~\invar.\invartheoremthree\ which states that
BA is distributed due to a phase separation invariant
measure if and only if  the second class particle paths belong to 
the set
$E:=\{e=(e_i)_{i\in \Z}\, :\, e_{2i}\in 2\Z\,\, \forall i\in \Z,{\rm\ and\ }
e_i-e_{i+1}\in \{-1,1\}\, \forall i\in \Z\}$
and have the distribution which is invariant with respect to the shift of $E$
by $2$ (i.e. $e_i\rightarrow e_{i+2}$). Above $e_i$ is interpreted as
the position of
the second class particle at time $i$.  For CA 184, a similar results holds
with the only difference that now the invariance is with respect to shift by 
$1$, and $E:=\{e=(e_i)_{i\in \Z}\, :\, e_i-e_{i+1}\in \{-1,1\}\, 
\forall i\in \Z$, and if for some $i$ and $j$, $j>i$ it holds that 
$e_i-e_{i-1}=-(e_{i+1}-e_i)=e_{i+2}-e_{i+1}=\ldots=
e_{j}-e_{j-1}=-(e_{j+1}-e_j)$ then $j-i$ is odd$\}$, 
other words, the second class  particle
motion in CA 184  is a time-homogeneous  process with the state space
$E$. 
We note that the ideas similar to those just presented, have been used in 
[\ferravi] to construct phase separation measures for a so-called 
Boghosian Levermore cellular automaton.

In Section~\hydrod\ we review and a generalize results 
from [\belfer]. 
Theorem~\hydrod.\hydrodtheoremone\ presents a hydrodynamic limit for the
shape of surface in SG. It says that if in SG $\{f_n(\cdot), n\in
\naturals\}$ 
the distribution of $c_nf_0(n\cdot)$
converges to some process $W(\cdot)$ for some sequence of numbers $\{c_n\}$
then $c_n f^\ast_n(n\cdot)$ converges to the process $W^{\min}(\cdot)$ which is
defined by $W^{\min}(x):=\min\{W(y)\, ,\, x-1\leq y\leq x+1\},\, \forall x\in
\reals$; above $f^\ast_n$ denotes  a particular function whose shape coincides
with the shape of $f_n$. The cornerstone of this result is the following 
property of SG: if for a function $f(\cdot)\in \cR$ we define the function
$g(\cdot)$ by $g(x)=\min\{f(y)\, ,\, x-1\leq y\leq x+1\}$ then $g(\cdot)$ and
$\hat f(\cdot)$ have the same shape. With help of 
Assertions~\introd.\introdassone, \introd.\introdasstwo\ we then derive from 
Theorem~\hydrod.\hydrodtheoremone\ the $time\rightarrow\infty$ 
 limits for the particle
distribution in BA (Theorem~\hydrod.\hydrodtheoremtwo) and in  CA 184
(Theorem~\hydrod.\hydrodtheoremthree), under appropriate rescaling
(i.e. hydrodynamic limits). These theorems express 
the limit of distributions  of particles in BA and CA~184
in terms of a limit of their counting
processes. A counting process for BA $\{\zeta_n, n\in \naturals\}$ is, 
by definition, the process $\{T_{BA, SG}(\zeta_n), n\in \naturals\}$. 
Since the latter closely relates to SG, as specified in
Assertion~\introd.\introdasstwo\ and Remark~\introd.\introdremarktwo, then
Theorem~\hydrod.\hydrodtheoremtwo\ is immediate from 
Theorem~\hydrod.\hydrodtheoremone. To obtain a counterpart for CA 184, just a
simple trick is required. The results of Section~\hydrod\  have two natural 
applications which
we demonstrate via examples. Example~\hydrod.\hydrodexampleone\ shows how BA 
and CA 184 may be employed in order to find the law of $W^{\min}(\cdot)$ for 
a given process $W(\cdot)$. Example~\hydrod.\hydrodexampletwo\ demonstrates a
simple way to 
find approximately the particle distribution in BA or CA 184 for large times.

In Section~\rate\ we calculate the rate of convergence of CA~184 starting 
from the Bernoulli 1/2 product measure, to its invariant state. This is 
equivalent to calculating the rate of the decay of particles in BA with 
a particular initial state. The argument we employ is not new. It has been 
used to estimate analogous rates in various
models (see [\cyclic], [\krugspohn], [\elskens]). This argument is based on
estimating a certain  quantity which we specify in Remark~\rate.\rateremarkone.
We observe (in Remark~\rate.\rateremarkzero) that the result discussed in this 
section could be also obtained using the tools presented in Section~\hydrod.

We consider the results of Theorems~\invar.\invartheoremthree,
\invar.\invartheoremfour, \hydrod.\hydrodtheoremthree\ to be our novel
contribution to the study of the dynamics of BA, CA~184 and
SG. Theorems~\invar.\invartheoremone, \invar.\invartheoremtwo\ and
\hydrod.\hydrodtheoremone, \hydrod.\hydrodtheoremtwo\ have original form but
the bulk of their proofs is based on ideas and techniques that have appeared
in literature related to the field.
\vskip4mm
\noindent{\bf \invar. INVARIANT MEASURES}\vskip2mm
\vskip2mm

In this section, we give the necessary and sufficient condition for a measure
to be invariant for BA (Theorem~\invar.\invartheoremone) and for CA 184
(Theorem~\invar.\invartheoremtwo); the corresponding result for the shape in 
SG may be
easily obtained from that for BA (by using $T$ from
Assertion~\introd.\introdasstwo)   and thus, will not be presented here.
Those invariant measures which are translation invariant are characterized in 
Remark~\invar.\invarremarktwo, and those which
are not translation invariant are characterized in 
Theorems~\invar.\invartheoremthree\ and \invar.\invartheoremfour.
At the end of Section~\invar\ we calculate the flux of
particles in CA 184 when it is distributed with an invariant and
translation invariant measure.

We start with some definitions. Elements from $\{0,1\}^{\Z}$ and from 
$\{-1,0,1\}^{\Z}$ will be called {\it configurations}, and their values at
sites of $\Z$ will be interpreted in terms of particles as indicated in
Introduction. A-particles of BA will be called simply {\it particles}, and
those whose velocity is $+1$ (resp., $-1$) will be called {\it positive}
(resp., {\it negative}). We set $\Theta:=\{-1,0,1\}^{\Z}$.  
By $\cI$ we will denote the set of the measures that are invariant
for the dynamics of BA. To say $\mu\in \cI$ means that
$ \mu\bigl(A^{-1}(U)\bigr)=\mu \bigl(U\bigr)$
for any cylinder set $U\subset \Theta$ (that is, a set of the form 
$U:=\{\zeta\in \Theta\, :\, \zeta(i_1)=a_1, \ldots, \zeta(i_k)=a_k\}$
for some $k\in \naturals$, $i_1, i_2, \ldots, i_k\in \Z$ and 
$a_1, a_2, \ldots, a_k\in \{-1,0,1\}$.) By $\tau_1$ we will denote the 
translation of $\Z$ by $1$ to the right. By $\cT$ we shall denote the set of
those measures on $\{-1,0,1\}^{\Z}$ that are translation invariant,
{\it i.e.},
$\mu(\tau^{-1}_1U)=\mu(U)$ for any cylinder set $U\subset \Theta$. 
By $\phi$ we shall denote the configuration that does not have any particle: 
$\phi(i)=0\, \forall i\in \Z$. By $\Theta_+:=\{0,1\}^{\Z}
\setminus\{\phi\}$ (resp., $\Theta_-:=\{-1,0\}^{\Z}\setminus\{\phi\}$) we
denote  the set of those  configurations  of 
$\Theta$ which have solely positive (respectively, negative)  particles.
By $\cP_+$ (resp., $\cP_-$) we denote  the set of measures concentrated on 
$\Theta_+$ (resp., on $\Theta_-$).
Given a configuration $\zeta\in \Theta$, we will call two particles of
$\zeta$ {\it subsequent}, if $\zeta$ does not have any other particle between 
these two. A pair of subsequent particles is said to be {\it converging}
if the leftmost particle of the pair is positive and the 
rightmost one is negative. The {\it diverging} type is introduced in 
analogous way. By $\Theta_s$ we denote the set of those configurations
$\zeta\in \Theta$ which satisfy the following two conditions
\vskip1mm
\item{{\it (i)}} there is only one pair
of subsequent converging particles;
\item{{\it (ii)}} both the number of positive particles and the number
of negative ones in $\zeta$ are infinite.
\vskip1mm
\noindent By $\cS$ we will denote the set of measures concentrated on 
$\Theta_s$. These measures are customarily called {\it phase separation 
measures}. 
\vskip2mm

\noindent{\bf Theorem~\invar.\invartheoremone.}\ \ {\sl $\mu \in \cI$ if and 
only if  for some  $\alpha_1, \alpha_2, \alpha_3, \alpha_4\geq 0$ such that 
$\alpha_1+\alpha_2+\alpha_3+\alpha_4=1$ it holds that
$$
\mu=\alpha_1\mu_1+\alpha_2\mu_2 +\alpha_3\mu_3+\alpha_4\delta_\phi, 
\eqno(\invar.\decomposition)
$$
whereas
$
\mu_1\in \cP_+\cap \cT, \,\,\, \mu_2\in \cP_-\cap \cT, \,\,\, \mu_3\in \cS\cap
\cI\cap\cT^c
$
.}
\vskip2mm

\noindent{\bf Remark~\invar.\invarremarkone.}\ \ $\cT^c$ means the complement
to  $\cT$. We wrote  $\cS\cap \cI\cap \cT^c$ in order to emphasize that  
$\cS\cap \cI$ cannot contain a translation invariant measure. This fact can be
established by the following reasoning: Let
$pos(\zeta)$ and  $neg(\zeta)$  denote the positions of respectively, the
positive and the negative  particles of the unique pair of converging
particles from  $\zeta\in \Theta_s$. Define the function 
$X:\Theta_s\rightarrow \Z\cup (\Z+ 1/2)$ by $X(\zeta):=
(neg(\zeta)-pos(\zeta))/2$. Now, if $\zeta$ is distributed by a translation
invariant measure then the distribution of $X(\zeta)$ would be invariant with
respect to translations of $\Z$, which is impossible.
\vskip2mm

\noindent{\bf Remark~\invar.\invarremarktwo.}\ \  
The set $(\cP_+\cap \cT)\cup \delta_\phi$ (and $(\cP_-\cap \cT)\cup 
\delta_\phi$, by analogy) is by its definition,
the set of all translation invariant measures on $\{0,1\}^{\Z}$. The
extremal points of this set is known to coincide with the set of all ergodic
(with respect to $\tau_1$) measures on $\{0,1\}^{\Z}$. 
\vskip2mm

\noindent {\it Proof of Theorem~\invar.\invartheoremone.}\ Observe that
$$
{\rm if\ }\zeta\in \Theta_+{\rm\  then\ }\tau_1(\zeta)=A(\zeta)
\eqno(\invar.\bum)
$$
It follows then easily from (\invar.\bum) that if $\mu_1\in \cP_+\cap \cT$ 
then $\mu_1\in \cI$.
Analogously, if $\mu_2\in \cP_-\cap \cT$ then $\mu_2\in \cI$. Since $\mu_3$
in (\invar.\decomposition) is picked from a subset of $\cI$ and since
obviously $\delta_\phi\in \cI$  then $\mu$  defined by (\invar.\decomposition)
must
belong  to $\cI$. This completes the proof of the ``if'' part of the theorem.
The ``only if'' part is more difficult. It will be proved with the aid of the
following 
\vskip1mm

\noindent{\bf Assertion~\invar.\invarassone.}\ \ 
{\sl Let $\zeta_0\in \Theta$ and $m\in \naturals$ be arbitrarily 
fixed. Let  $\zeta_m$ denote the configuration of particles in the BA
at time $m$, starting from $\zeta_0$, {\it i.e.}, $\zeta_m=A^m\zeta_0$. 
Then the distance between the particles of any pair of diverging particles
from $\zeta_m$ is not less than $2m+1$.} 
\vskip2mm

\noindent{\it Proof of Assertion~\invar.\invarassone.} \ 
Find any pair of diverging  particles in $\zeta_m$. Observe that the initial
position of the negative particle of this pair must be to the left of that of
the positive one, because if not then these particle would have annihilated
each other by time $m$. The assertion follows from this observation and from 
the fact that these particles diverge with velocity $2$.
\hfill$\clubsuit$
\vskip2mm

\noindent{\it Continuation of the proof of Theorem~\invar.\invartheoremone}. \ 
 Assertion~\invar.\invarassone\ provides that
if $\mu\in \cI$ then it gives weight zero to any configuration that has at
least one pair of diverging particles. Thus such $\mu$ must be
concentrated on $\Theta_+\cup \Theta_-\cup \Theta_s\cup\{\phi\}$. Since the
dynamics $A$ does not mix these sets (that is, there is no a configuration
$\zeta$ from one of these sets such that $A\zeta$ belongs to another one)
then the expansion (\invar.\decomposition) follows. \hfill$\clubsuit$
\vskip2mm

Let $T$ be defined by (\introd.\relation).
We define $\Lambda\subset \Theta$ as the image of $\{0,1\}^{\Z}$ by
$T$. Namely, $\Lambda:=\{\zeta\in \Theta\, :\, \exists  \eta\in
\{0,1\}^{\Z}$ such that $T\eta=\zeta\}$. 
>From this definition and that of $T$,
$\Lambda$ is easily characterized: $\Lambda:=\{\zeta\in \Theta:$ the distance
between any subsequent particles with the same velocity is odd, and the
distance between any two subsequent particles with opposite velocities is
even$\}$. The set $\Lambda$ is used in the following statement:
\vskip2mm

\noindent {\bf Theorem \invar.\invartheoremtwo.}\ \ 
{\sl A measure $\mu$ on $\{0,1\}^\Z$ is invariant for CA~184 
if and only if 
$$
\mu=\alpha_1\mu_1+\alpha_2\mu_2 +\alpha_3\mu_3+\alpha_4\biggl({{\delta_o+
\delta_e}\over{2}}\biggr), 
{\sl\ for\ some\ } \alpha_1, \alpha_2, \alpha_3, \alpha_4\geq 0, \,
\alpha_1+\alpha_2+\alpha_3+\alpha_4=1
$$
where the measures  $\mu_1$ and $\mu_2$ are translation invariant and are 
concentrated 
on the sets $T^{-1}(\Theta_+\cap \Lambda)$ and $T^{-1}(\Theta_-\cap \Lambda)$
respectively; the measure $\mu_3$ is invariant for CA~184 but
not translation invariant and is concentrated on $T^{-1}(\Theta_s\cap 
\Lambda)$; the measure $\delta_o$ (resp., $\delta_e$) gives mass $1$ to the
configuration $o$  (resp., $e$) where $o$ and $e\in \{0,1\}^{\Z}$ are such
that $o(i)=1,\, e(i)=0$, when $i$ is odd, and $o(i)=0,\, e(i)=1$ when $i$ is
even, $i\in \Z$.
}
\vskip2mm

\noindent{\it The proof of Theorem~\invar.\invartheoremtwo}\ comes out by a
reformulation of the proof of Theorem~\invar.\invartheoremone\ in terms of
CA~184. This reformulation is based on the property of $T$ stated in
Assertion~\introd.\introdassone, and the
facts that $T\bigl(\{0,1\}^{\Z}\bigr)=\Lambda$ and that $T$ is invertible on
$\{0,1\}^{\Z}\setminus \{o,e\}$. The only care thus, should be taken because
of the fact that $T(o)=T(e)=\phi$. But this causes no difficulties in the
mentioned reformulation because we have that $C(o)=e$ and $C(e)=o$ and as a
consequence of this, we also have that if $\mu$ is invariant for CA~184
then $\mu(o)=\mu(e)$.
{}\hfill$\clubsuit$
\vskip2mm

We  now characterize the set of invariant measures of BA which are not
translation invariant. Set 
$$
E:=\{e=(e_i)_{i\in \ZZ}\, :\, e_i\in \Z\,\, \forall i\in \Z{\rm\ and\ }
e_i-e_{i+\half}\in \{-\half, \half\}\,\, \forall
\, i\in \ZZ\}
$$
and define the shift $\tau_E\,:\, E\rightarrow E$  by
$\tau_E(e)=e^\prime\, \Leftrightarrow\,e^\prime_i=e_{i+1}\, \forall i$. 
Denote then by
$\cE_\tau$ the set of the measures on $E$ which are invariant with respect to
the shift $\tau_E$.
\vskip2mm

\noindent{\bf Theorem~\invar.\invartheoremthree.}\ \ {\sl 
There is a bijection $F^\ast:\cS\cap \cI\rightarrow \cE_\tau$.
The mapping $F^\ast$ will be constructed explicitly in the proof of the
theorem.}
\vskip2mm

\noindent{\it Proof.}\ \ The idea is to introduce a second class
particle in BA and
to show that BA is distributed by an invariant phase separation measure if and
only if the second class particle's motion has invariant distribution with 
respect to
shifts of time. The  second class particle in
BA is defined by the following rules: {\it (i)}\ it is an extra particle
that does not affect the evolution of other particles, and {\it (ii)}\ it
moves with the velocity that changes from $+1$ to  $-1$ and back depending on
the BA particles that the second class particle meets, namely:{\it (ii-a)}\
when the second class
particle meets a negative particle it changes its velocity to $-1$ and keeps 
going together with the met particle; {\it (ii-b)}\ at the time this met 
particle is annihilated, the second class particle changes its velocity 
to $+1$. 
One more postulate is: {\it (iii)}\ if two particles annihilate at a time 
$n\in \ZZ$ then at this time, the second class particle is 
at the annihilation point and changes its velocity from $-1$ to $+1$.
We want the  rules {\it (i)-(iii)}\ to define correctly and uniquely the
second class
particle and moreover, we want each trajectory of the second class particle
 to determine
uniquely the configuration of BA at any time. To achieve these objectives we 
have to consider BA with double-infinite time and to record the second class
particle
position at times from the set $\ZZ$. Note that $\Z+\half$ has to be added
because annihilations may occur at times from this set. Also note that the
distribution of the second class particle trajectories turns out to be 
invariant 
with respect to
the shift of time by $1$ but not necessarily to the shift of time by $\half$.
The details are presented below. 

We start introducing a terminology which we shall need to construct
BA with double infinite time. For an element $\o$ from
$\{-1,0,1\}^{\Z\times\Z}$, let us write $\o_n(m)$ for its values at $(n,m)\in
\Z\times \Z$. Let us also write $\o_n$ for the restriction of $\o\in
\{-1,0,1\}^{\Z\times\Z}$ obtained by fixing its first coordinate to $n$;
observe $\o_n\in \{-1,0,1\}^{\Z}$. Define then $\O:=\{\o\in
\{-1,0,1\}^{\Z\times \Z}\,:\, \o_n\in \Theta_s{\rm\ and\ }
 \o_{n+1}=A\o_n\, \forall n\}$. Fix now 
an arbitrary $\mu\in \cS\cap \cI$. There is a standard procedure to construct
a probability measure $\pr^\mu_\O$ (on an appropriate
$\sigma$-algebra $\cF$ of subsets of $\Omega$) with the following properties
(\invar.\propone), (\invar.\proptwo):
$$
\pr^\mu_\O\bl \o\, :\, \o_n\in D \br=\mu\bl D \br, 
\kern2em\forall n\in \Z{\rm\ and\ } 
\forall {\rm\ cylinder\
} D\subset
\{-1,0,1\}^{\Z}
\eqno(\invar.\propone)
$$
\noindent If we define the shift operator $\tau_\O:\Omega\rightarrow\Omega$ 
acting by $(\tau_\O(\o))_n(m):=\o_{n-1}(m), \forall n,m$, then 
$$
\pr^\mu_\O \bl D\br=\pr^\mu_\O\bl \tau_\O^{-1}(D)\br,\kern2em \forall 
{\rm\ cylinder\ } D\subset \{-1,0,1\}^{\Z\times \Z}
\eqno(\invar.\proptwo)
$$
$\pr^\mu_\O$ will be called {\it BA with double infinite time and the marginal
distributions $\mu$}. Such a choice of the name becomes obvious if one 
interprets $\o$ as a trajectory of a certain process, and $\o_n(m)$ as its 
value at time $n$ and the location $m$. We recall that it is known that given
$\mu\in \cS\cap \cI$, there is a unique $\pr_\O^\mu$ satisfying
(\invar.\propone),  (\invar.\proptwo), and,
moreover, given any probability measure $\pr$ on $\O$ which is invariant with
respect to $\tau_\O$, there is a unique $\mu\in \cS\cap \cI$ such that
$\pr_\O^\mu=\pr$, other words, $\mu\rightarrow \pr_\O^\mu$ is a bijection from
$\cS\cap\cI$ to the set of all probability measures on $\O$ which are
 invariant in respect to $\tau_\O$. 

By certain technical reasons that will become clear soon, we need to have 
a record 
of the positions of  particles in BA also at times from the set
$\Z+\half$. We thus associate to each $\o\in \O$ an element $\o^\prime\in
\{-1,0,1\}^{\ZZ\times\ZZ}$ by the following rule: for each $n\in \Z$ we set 
$$
\o^\prime_n(m)=\cases{\o_n(m),& when $m\in \Z$ \cr
                        0, & when $m\in \Z+\half$\cr
}\eqno(\invar.\primeone)
$$
and
$$
\o^\prime_{n+\half}(m)=\cases{ 0, & if $m\in \Z$\cr
 1, & when $m\in \Z+\half$ and if $\o_n(m-\half)=1$ and
$\o_n(m+\half)\not=-1$\cr 
-1, & when $m\in \Z+\half$ and if $\o_n(m+\half)=-1$ and 
$\o_n(m-\half)\not= 1$\cr
                    0, & otherwise\cr}\eqno(\invar.\primetwo)
$$
The relation between $\o$ and $\o^\prime$ may be explained in the following
manner: Suppose that two persons observe the evolution of the same process
BA. Suppose that one of them records the states of the process at times 
$\Z$ while the other one records the states of the process at times 
$\ZZ$. Then, $\o$ and $\o^\prime$ are the histories recorded by respectively,
the first and the second person. Having this interpretation in mind and
observing that a configuration of particles in BA at time $n\in \ZZ$ 
determines uniquely the configuration at time $n+\half$ we conclude that 
the mapping $\o\rightarrow\o^\prime$ is a bijection.

We are now in a position to define the mapping $F:\O\rightarrow E$. 
To obtain $e=F(\o)$, we
construct first $\o^\prime$ by the rules (\invar.\primeone)-(\invar.\primetwo)
and then,
for each $n\in \ZZ$, we define $e_n$ by the following rules:

\item{1} if a point $m\in \ZZ$ is an annihilation place of two particles at
time $n$ in $\o^\prime$ (that is, if $\o^\prime_n(m)=0,
\o^\prime_{n-\half}(m-\half)=1, \o^\prime_{n-\half}(m+\half)=-1$) 
then we set $e_n=m$;

\item{2} if $e_{n-\half}=m$ and $\o^\prime_{n-\half}(m)=-1$ then we set 
$e_{n}=e_{n-\half}-\half$;

\item{3} if $e_{n-\half}=m$ and $\o^\prime_{n-\half}(m)=0$ then we set
 $e_{n}=e_{n-\half}+\half$.
\vskip1mm

\noindent An important observation is that $F$ is a bijection. This follows
from two facts. The first one is that $\omega\rightarrow\o^\prime$ is a 
bijection, as mentioned above. The second fact is that $\o^\prime$ determines
uniquely $e\in E$ and vice versa. This fact may be verified straightforward. 
We observe that in order that $F^{-1}e$ be uniquely determined, it was 
essential  that the elements of $E$ be indexed by the set $\ZZ$. (We suggest a
reader to choose an arbitrary $\epsilon\in E$, to construct 
then $\varepsilon\in E$
by changing the value of $\epsilon$ at $\half$, and to verify then that 
$F^{-1}(\epsilon)\not=F^{-1}(\varepsilon)$.) 
The mapping $F$ introduces naturally a measure $\pr^\mu_E$ on $E$ via
$\pr^\mu_E\bl D\br:=\pr^\mu_\O\bl F^{-1}D\br, \forall$ cylinder $D\subset E$. 
But because of the property
$$
\tau_\O\bigl(F^{-1}(e)\bigr)=F^{-1}\bigr(\tau_E (e)\bigr)
$$
that follows from our constructions and because of (\invar.\proptwo), we
conclude that $\pr^\mu_E\bl D\br=\pr^\mu_E\bl \tau^{-1}_E D\br, 
\forall D\in E$, 
that is $\pr^\mu_E$ is from $\cE_\tau$. We thus define $F^\ast$ by
$F^\ast(\mu):=\pr^\mu_E$, for $\mu\in \cS\cap\cI$. The fact that $F^\ast$ is a
bijection follows from that the mappings $F$ and $\mu\rightarrow\pr^\mu_\O$
are bijections.\hfill$\clubsuit$
\vskip4mm

Recall Theorem~\invar.\invartheoremtwo. Let us give the name {\it invariant 
phase
separation measure of CA 184} to an invariant measure of CA 184 which is
concentrated on $T^{-1}(\Theta_s\cap \Lambda)$. The following result
characterizes these measures.
\vskip2mm

\noindent{\bf Theorem~\invar.\invartheoremfour.}\ \ {\sl There is a one-to-one
mapping between the set of all invariant phase separation measures of CA
184 and the  set of all stationary processes $(X_n)_{n\in \Z}$ with values 
in $\Z$ which satisfy the conditions
\item{(a)} $X_{n+1}-X_n\in \{-1,1\}$, for all $n$;
\item{(b)} if for some $i, j$, $j>i$, it holds that
$X_i-X_{i-1}=-(X_{i+1}-X_i)=X_{i+2}-X_{i+1}=\ldots=
X_{j}-X_{j-1}=-(X_{j+1}-X_j)$ then $j-i$ is odd.

\noindent This mapping can be obtained as an appropriate  restriction of 
$F^\ast$ in the
manner outlined in the proof.}
\vskip2mm

\noindent{\it Proof.}\ \  The mapping $T$ from
(\introd.\relation) induces a natural mapping between the distributions of
CA~184 and that of BA. The property 
that $T:\{0,1\}^{\Z}\rightarrow \Lambda$ is a
bijection (recall, $\Lambda$ has been defined before the statement of
Theorem~\invar.\invartheoremtwo) then yields that the induced mapping is a
bijection between  the set of all invariant phase separation measures of 
CA~184 and a subset of $(\cS\cap\cI)$ which we shall denote by $(\cS\cap
\cI)_{{\rm 184}}$. 

Thus, Theorem~\invar.\invartheoremfour\ is proved, whence we show that 
$F^\ast\bigl((\cS\cap\cI)_{{\rm 184}}\bigr)$ is the set of the processes $X$
as defined in the theorem's statement. This can be done directly via repeating
the proof of Theorem~\invar.\invartheoremthree, but applying it to 
$(\cS\cap\cI)_{{\rm 184}}$ instead of $\cS\cap \cI$. For this application to
be correct, it is necessary to note the difference between these two sets.
It follows from the property of $T$ mentioned in the above paragraph that
a measure from
$\cS\cap \cI$ is in $(\cS\cap \cI)_{{\rm 184}}$ if and only if it gives weight
$0$ to the set $\Lambda^c=\{\zeta\in \{-1,0,1\}^{\Z}\,\,$: the distance
between any subsequent particles with the same velocity is even or  the
distance between any two subsequent particles with opposite velocities is
odd$\}$. This suggests the following two modifications in construction of 
$F^\ast\bigl((\cS\cap\cI)_{{\rm 184}}\bigr)$ compared to that of 
$F^\ast\bigl(\cS\cap\cI\bigr)$: First, we do not need to record the second
class particle positions at times $\Z+\half$, because there are no
annihilations at these times. Second, the set of  paths of the second class 
particle should be modified so that when a path from this set 
is translated into corresponding configuration of annihilation
particles, then the latter must belong to $\Lambda^c$. To satisfy these
requirements, we imposed the conditions {\sl (a) and (b)} on the process $X$. 
The rest of the proof goes almost identically to that of
Theorem~\invar.\invartheoremthree\ and thus, will be omitted.\hfill$\clubsuit$
\vskip4mm

\noindent {\bf Remark~\invar.\invarremarkfour.} \  We recall 
(from [\liggett] and [\japonez],
respectively) that in the case of TASEP and
CA\&N, the set of phase separation measures is $\{ \nu^{(n)},
-\infty< n<\infty\}$, where $\nu^{(n)}$ gives mass $1$ to the configuration
$\eta^{(n)}$ such that $\eta^{(n)}(x)=1\,\forall\, x\geq n$ and
$\eta^{(n)}(x)=0\,\forall\, x< n$.
\vskip2mm

\noindent{\bf Flux of particles in CA 184} which is distributed due to an
translation invariant measure $\mu$ which is invariant for CA~184,
is defined as
$$
\mu\{\eta_n(1)(1-\eta_{n-1}(1))\}
$$
We shall now show that the flux of particles in CA~184 is given by 
$1/2-|1/2-\rho|$, if CA~184 distributed by
such a measure $\mu$, where $\rho$ denotes the density of particles.
(Here the density of particles is understood as
$\mu\{\eta_n(i)\}$ which certainly does not depend on $n$ and $i$.) Let $\mu$
be such a measure and consider first the case $\rho<1/2$. Then
Theorem~\invar.\invartheoremtwo\ implies that $\mu\in T^{-1}(\Theta_-\cap
\Lambda)$ which in turn, implies that $\mu\{\eta_n(i)=0\given \eta_{n}(i-1)=1\}
=1$, $\forall i\,\, \forall n$. Thus the flux of particles is
$\rho$. In the case when $\rho>1/2$, the flux of particles is $1-\rho$. This
results from the following reasoning: We define the flux
of holes as $\mu\{(1-\eta_n(0))\eta_{n-1}(0)\}$ and derive that it equals
$1-\rho$, this result follows by adapting the reasoning that have been used
for the case $\rho<1/2$. We then observe that the flux of holes must be equal
to the flux of particles (this observation may be justified in one of the
following ways: either using the fact that $\mu$ is invariant, or from a direct
comparison of the definitions of the fluxes). It is left to consider the case
$\rho=1/2$ which occurs as follows from Theorem~\invar.\invartheoremtwo, when
$\mu=1/2(\delta_o+\delta_e)$. When it happens, then flux is obviously $1/2$.
\vskip4mm


\noindent{\bf \hydrod. HYDRODYNAMIC LIMIT.}
\vskip2mm

This section concerns the limit (as $time\,\rightarrow\infty$) shape of surface
in SG (Theorem~\hydrod.\hydrodtheoremone) and the limit distribution of
particles in BA (Theorem~\hydrod.\hydrodtheoremtwo) and in CA 184
(Theorem~\hydrod.\hydrodthethree). These results are called hydrodynamic
limits because  along with $time$ tending to $\infty$ we apply a certain
rescaling that depends on $time$, to the quantities whose temporal limit is of
interest. We also demonstrate here applications
of these results (Examples~\hydrod.\hydrodexampleone,
\hydrod.\hydrodexampletwo). 

For  $y\geq 0$, let us denote by  $M_y$ the operator that brings a
function $f(\cdot):\reals\rightarrow\reals$ to the function
$g(\cdot):\reals\rightarrow\reals$ in the way such that
$$
g(x)=\bigl(M_yf\bigr)(x)=\min\{f(z)\, :\, x-y \leq z\leq x+y\}, \kern2em
\forall x\in \reals
$$
Let us then introduce {\it Surface Modification Process} (SM) as the name of a
sequence $\{f_n(\cdot), n\in \naturals\}$ of random functions from $\cR$ such
that $f_{n+1}(\cdot)=(M_1f_n)(\cdot)$ for each $n\in \naturals$. One then can
check straightforwardly the following property:
$$
{\rm for\ any\ }f\in \cR {\rm\ there\ is\ an\ }\alpha\in \reals{\rm\  such\ 
that\ }(M_1f)(\cdot)=\hat f(\cdot)+\alpha\eqno(\hydrod.\efen) 
$$
whereas the construction rule of $\hat f(\cdot)$ has been specified in the 
part of Introduction that defined the Surface growth process.
The above (\hydrod.\efen) says that at each time, the shape of the states 
of SG and SM 
is the same, provided they stated from the same state. However, SM is more
handy for formulating the hydrodynamic limit for this shape. 
To state this result, we need first to introduce some notations: {\it
(1)}\ For a process $W(x), \, x\in \reals$, we will denote by 
$W^{\min}(\cdot)$ the
process defined by $W^{\min}(x):=\bigl(M_1 W\bigr)(x), \,\forall x\in \reals$.
{\it (2)}\ 
The convergence of processes discussed in this section and throughout  the
paper,  is understood as the weak convergence of their distributions on
$\cC[a,b]$, the set of continuous functions on $[a,b]$, for any
$-\infty<a<b<\infty$ (see [\bilin]).
\vskip2mm

\noindent{\bf Theorem~\hydrod.\hydrodtheoremone}\ (hydrodynamic limit for the
shape of surface in SG). \
{\sl Let $f_0(\cdot)$ be a random function with the state space $\cR$ such that
$c_nf_0(n\cdot)\rightarrow W(\cdot)$, as $n\rightarrow\infty$, for some
sequence $\{c_n,\, n\in \naturals\}$ of real numbers and some process $W(x),
x\in \reals$. Then the SM $\{f_n(\cdot),\, n\in \naturals\}$ 
satisfies the following property 
$$
c_nf_n(n\cdot)\rightarrow W^{\min}(\cdot)\eqno(\hydrod.\state)
$$ 
}
\vskip2mm

\noindent{\it Proof.}\ \ First, we note that for any function $h:\reals
\rightarrow\reals$, the modulus of continuity (which definition one finds in
[\bilin]) of $M_1h$ does not exceed that
of $h$; this property of $M_1$ can be verified straightforward. This fact
allows us to apply  Theorem~5.1 from [\bilin] to the assumption $c_nf_0(n\cdot)
\rightarrow W(\cdot)$ to conclude that 
$$
M_1\bigl(c_nf_0(n\cdot)\bigr)(\cdot)\rightarrow  \bigl(M_1 W\bigr)(\cdot)
=W^{\min}(\cdot) \eqno(\hydrod.\frombilin) 
$$
We must comment on the usage of two "$\cdot$" in (\hydrod.\frombilin) in order
to avoid a possible confusion in its interpretation: $M_1(c_nf_0(n\cdot))(\cdot)$
means the function obtained  by rescaling $f_0$ by $c_n$ along the axis of 
ordinates and by $n^{-1}$ along the axis of abscissas, and then, by applying  
$M_1$.
It may be verified directly that for any continuous function $h$ and
any constant $c$,
$$
\bigl(M_1 c h(n\cdot)\bigr)(x)=c\bigl(M_n h \bigr)(nx), \kern2em
\forall x\in \reals,
$$ 
and also that $n$ subsequent applications of $M_1$ is equivalent to one 
application of $M_n$. Thus, (\hydrod.\state) follows from
(\hydrod.\frombilin). \hfill$\clubsuit$
\vskip4mm

Let $(\Omega, {\cal F}, \pr)$ be an abstract probability space. 
We shall say that a random function
$f_n(\cdot):\Omega\rightarrow \cR$ is {\it a counting profile} of
$\zeta_n:\Omega\rightarrow \{-1,0,1\}^{\Z}$ (resp., $\eta_n:\Omega \rightarrow
\{0,1\}^{\Z}$) if $f(k)[\o]-f(k-1)[\o]=\zeta(k)[\o]$ (resp.,
$f(k)[\o]-f(k-1)[\o]=-\eta(k)[\o]$), $\forall k\in \Z$ (certainly,
it is enough to specify the values of a function from $\cR$ just on $\Z$) and 
$\forall \o\in \Omega$, where we denoted by $Y[\o]$ the realization on 
$\o\in \Omega$ 
of a function $Y$ defined on $\Omega$. We say that $\{f_n(\cdot), n\in 
\naturals\}$ is {\it a counting process} of BA $\{\zeta_n, n\in \naturals\}$ 
(resp., CA 184 $\{\eta_n, n\in \naturals\}$) if for each $n$, $f_n(\cdot)$ is 
a counting profile of $\zeta_n$ (resp., $\eta_n$). 

>From the definitions just introduced, Assertion~\introd.\introdasstwo\ and
(\hydrod.\efen), we conclude that if a random function
$f_0(\cdot)$ is a counting profile of a random configuration $\zeta_0$ 
then SM $\{f_n(\cdot), n\in \naturals\}$ is a counting 
process for  BA $\{\zeta_n, n\in \naturals\}$. Thus,
Theorem~\hydrod.\hydrodtheoremone\ leads to the following
\vskip2mm
 
\noindent{\bf Theorem~\hydrod.\hydrodtheoremtwo}\ (hydrodynamic limit for
particle distribution in BA). \ 
{\sl Let $\{\zeta_n, n\in \naturals\}$ be a BA. Assume there is $f_0(\cdot)$,
a counting profile of $\zeta_0$, such that $c_nf_0(n\cdot)\rightarrow
W(\cdot)$ as $n\rightarrow\infty$, for some sequence $\{c_n,\, n\in
\naturals\}$  of real numbers and some process $W(x),
x\in \reals$. Then there is a counting process $\{f_n(\cdot), n\in
 \naturals\}$ of BA  such that $c_nf_n(n\cdot)\rightarrow
W^{\min}(\cdot)$. Precisely to say, this counting process is SM with the
initial state $f_0(\cdot)$.}
\vskip2mm

The next result is a counterpart of Theorem~\hydrod.\hydrodtheoremtwo\ for
CA~184.  Note  that it contains  an additional assumption ``$c_n\rightarrow 0$
as $n\rightarrow\infty$''. In respect to this note, we observe that in all 
applications of hydrodynamic type limit theorems we are aware of, this 
condition is satisfied.
\vskip2mm

\noindent{\bf Theorem~\hydrod.\hydrodthethree}\ (hydrodynamic limit for
particle distribution in CA 184). \ 
{\sl Let $\{\eta_n, n\in \naturals\}$ be a CA 184. Assume there is 
$g_0(\cdot)$,
a counting profile of $\eta_0$, such that $c_ng_0(n\cdot)\rightarrow
W(\cdot)$ as $n\rightarrow\infty$, for some process $W(x),
x\in \reals$, and some sequence $\{c_n,\, n\in
\naturals\}$  of real numbers satisfying $c_n\rightarrow 0$ as $n\rightarrow
\infty$. Then there is a counting process $\{g_n(\cdot), n\in
\naturals\}$ of CA 184 such that $c_ng_n(n\cdot)\rightarrow
W^{\min}(\cdot)$. The construction rule of $g_n(\cdot)$ will be specified in
the proof of this theorem.}
\vskip2mm

\noindent{\it Proof.}\ \ It is straightforward to see that
$$\eqalign{
&{\rm if\ }T\eta=\zeta{\rm\ for\ some\ }\eta\in \{0,1\}^{\Z}{\rm\ and\ }\zeta
\in \{-1,0,1\}^{\Z},{\rm\ where\ }T{\rm\ is\ from\ (\introd.\relation),}\cr
&{\rm and \ if\ }g{\rm\ and\ }f{\rm\ are\ counting\
profiles\ of\ }\eta{\rm\ and\ }\zeta{\rm\ respectively,}\cr
&{\rm satisfying\ }g(0)=f(0){\rm\ then\ } |f(x)-g(x)|\leq 1\, \forall x\in
\reals\cr
}\eqno(\hydrod.\closeness)
$$

Let $\zeta_0:=T\eta_0$, where $T$ is from (\introd.\relation), and let
$f_0(\cdot)$ be the counting profile of $\zeta_0$ such that
$f_0(0)[\o]=g_0(0)[\o], \forall \o$. By theorem's assumptions and
(\hydrod.\closeness), $c_nf_0(n\cdot)\rightarrow W(\cdot)$. The latter and
Theorem~\hydrod.\hydrodtheoremtwo\ yield that  $c_n f_n(n\cdot)\rightarrow
W^{\min}(\cdot)$ for an appropriate counting process $\{f_n(\cdot), n\in
\naturals\}$ of the BA $\{\zeta_n, n\in \naturals\}$.
This conclusion, (\hydrod.\closeness) and the assumption $c_n\rightarrow 0$ then
yield that $c_ng_n(n\cdot)\rightarrow W^{\min}(\cdot)$, if we define 
$g_n(\cdot)$ to be the counting profile of $\eta_n$
such that $g_n(0)[\o]=f_n(0)[\o],\, \forall \o$. \hfill$\clubsuit$
\vskip2mm

\noindent{\bf Two applications of the above results} are quite natural. We
shall demonstrate them in examples below.
\vskip2mm

\noindent{\bf Example~\hydrod.\hydrodexampleone} demonstrates that given 
a process $W(\cdot)$,  the
distribution of $W^{\min}(\cdot)$ can be found by considering an appropriate
cellular automaton. This example is a brief review of a result contained in
[\belfer] in which $W(\cdot)$ was taken to be $B(\cdot)$, the
one-dimensional Brownian motion with $B(0)=0$.

Let $\zeta_n, n\in\naturals$, be the BA whose initial
state has the following distribution: 
$$
\zeta_0(i), i\in \Z, {\rm\ are\ i.i.d.\ with\ }\pr\bl \zeta_0(i)=1\br
=\pr\bl \zeta_0(i)=-1\br=1/2\eqno(\hydrod.\bernoulli)
$$
It is then not hard to find exactly the distribution of
$\zeta_n=A^n\zeta_0$. It is as follows:  it is invariant with respect to 
translations and reflections of $\Z$, and given a site
$i$ is occupied by a particle, the probability that the nearest particle to
its right has the same (opposite) velocity is $(1+u_{2n})^{-1}$ (resp.,
$u_{2n}(1+u_{2n})^{-1}$) and can be expressed exactly as a function of
$u_{2n}$, where $u_{2n}$ is the probability that a simple symmetric one
dimensional random walk returns to the origin at time $n$; this expression will
depend on the velocities of the particles under the consideration. For details
we refer a reader to Theorem~1 in [\belfer]. Two fact contributed to the
success in finding the law of $\zeta_n$ in such an exact form. First,
a simple expression for the distance between a given particle and its
annihilating companion (this is the essence of
Assertion~\rate.\rateassertionthree\  
presented in Section~\rate). Second, the relation of the distribution of this
distance to $u_{2n}$ which stems from a particular choice of the law of
$\zeta_0$, (\hydrod.\bernoulli).

Let now $f_n(\cdot)$ denote the counting profile of $\zeta_n$ that satisfies
$f_n(0)=0$. Clearly, the law of $f_n(\cdot)$ can be derived from that
of $\zeta_n$ presented above. It turns out that the known asymptotics for 
$u_{2n}$ (see [\feller]) is sufficiently precise to enable one to derive then
the distribution of the process to which $n^{-1/2}f_n(n\cdot)$ converges
(for details, see Theorem~2 in [\belfer]). If $f_n$ were $M_nf_0$ then this
law would be of $B^{\min}(\cdot)$, by Theorem~\hydrod.\hydrodtheoremtwo\
and because  $n^{-1/2}f_0(n\cdot)\rightarrow B(\cdot)$. By technical reasons 
however (recall Remark~\introd.\introdremarktwo) it is not. Nevertheless, this
limit law allows to expose the sample paths of $B^{\min}(\cdot)$. They are
described below following [\belfer].

A generic trajectory of  the process $\bm (\cdot)$ is a continuous function
and it is composed of decreasing,
increasing and flat portions which follow each other in the following cyclic
order:  flat-increasing-flat-decreasing. A flat portion preceded
by a decreasing portion and followed by an increasing one has length $1$; 
such a portion is called {\it valley}. A flat portion preceded by an increasing
portion and followed by a decreasing portion, has the length $\hat\theta$ 
picked from the distribution  $\pr \bl \hat \theta\leq x\br=2x^{1/2}(1+x)^{-1}, \,
x\in [0,1]$. Such portion is called {\it a plateau}. An increasing portion is
constructed from a realization of an increasing process ${{\cal L}}(x), 
x\geq 0$, which
is stopped at the height that has an exponential-1 law. The process 
${{\cal L}}(\cdot)$
is a subordinator whose L\'evy  measure $\mu$ (see [\karatzas], Chapter~6) 
has the following form:
$$
\mu(d\ell)={{d\ell}\over {2(\pi \ell^3)^{1/2}}}, \,\, \ell\in (0,1), \,
d\ell\subset [0,1]
$$
A decreasing portion being taken with the sign ``-'', is a probabilistic 
replica of an increasing
portion. All the portions that compose $\bm(\cdot)$ are independent which
means that the random variables involved in the construction of these
portions should be taken to be independent.

We observe that what is lacking to have a  complete description of $\bm(\cdot)$
is the dependence between $\bm(0)$, the value of this process at  $x=0$, and
the form of $\bm(\cdot)$  as a whole. This happened because $\zeta_n$
preserves only the form of $(M_n f_0)(\cdot)$. We may present results
in this respect in the future ([\moving]). 
\vskip2mm

\noindent{\bf Example~\hydrod.\hydrodexampletwo.}\ \  We shall use here 
our knowledge of the law  of
$\bm(\cdot)$ in order to find approximately the distribution at time $n$ of
CA~184 with a particular initial distribution.

Let $\eta_n, n\in \naturals$ be a CA 184 whose
initial state has the Bernoulli 1/2 distribution, that is
$$
\eta_0(i), i\in \Z, {\rm\ are\ i.i.d.\ with\ }\pr\bl \eta_0(i)=1\br=\pr\bl
\eta_0(i)=0\br=1/2\eqno(\hydrod.\bernoulliautomat)
$$
Let next $g_0(\cdot)$ denote the counting profile of $\eta_0$ such that
$g_0(0)\equiv 0$. Then as it is well known, 
$
n^{-1/2}g_0(n\cdot)\rightarrow B(\cdot)
$,
where $B(\cdot)$ is as specified in Example~\hydrod.1 above. Thus, by
Theorem~\hydrod.\hydrodthethree,
$$
n^{-1/2}g_n(n\cdot)\rightarrow B^{\min}(\cdot)\eqno(\hydrod.\inapplication)
$$
for an appropriate sequence of counting profiles.

Our programme is now as follows:
We shall base on (\hydrod.\inapplication) to assume that
$$
g_n(\cdot)=n^{1/2}\bm(n^{-1}\cdot){\rm\ in\ distribution}\eqno(\hydrod.\true) 
$$
and shall combine this assumption with the law of $\bm(\cdot)$ presented in
Example~\hydrod.1 above, in order to derive the picture of the particle distribution in
$\eta_n$ for large $n$.

Although $g_n(\cdot)$ is never a constant, $n^{-1/2}g_n(\cdot)$ looks like a
constant on an interval $I\subset \reals$ in which either each odd site of
$\eta_n$ is occupied by a particle while each even site is empty, or vice
versa.  Let us say that in this case we see a duce
pattern of $\eta_n$ in $I$. Similar considerations suggest us to call
$g_n(\cdot)$ increasing  (resp., decreasing) on an interval
$I\subset \reals$ if  $\eta_n$  does not have a pair of contiguous
sites of $I$ occupied by particles (resp., holes) and at both endpoints one
finds in $\eta_n$ a pair of contiguous sites occupied by holes (resp., particles)
(the term ``hole'' used here, is a natural nickname for
site free of particle). Let us say that in this case we see a hole dominated
(resp., particle dominated) pattern of $\eta_n$ in $I$

Being equipped with the terminology of the above paragraph, we can now present
the result of the  programme we have layed down. It is the following picture of
the distribution of $\eta_n$:
One will see the space $\reals$ as divided in intervals in the way such that
each one of them is either a particle dominated, or a hole dominated, or a
duce pattern of $\eta_n$.  The interval endpoints are distributed as
homogeneous  point process on $\reals$ (whose parameters can be calculated
explicitly from the law of $\bm(\cdot)$). The intervals follow
each other in the cyclic order: particle dominated -- duce --
hole dominated -- duce --  particle dominated. The mean extension 
of an interval is $Cn +o(n)$, where $C=C_1$ for both a particle dominated and a
hole dominated regions, $C=C_2$ for a duce region following a particle
dominated region, and $C=C_3$ for a duce region preceding a particle
dominated region.  The mean number of the pairs of contiguous sites occupied
by particles (resp., free of particles) in a particle dominated (resp., hole
dominated) region is $C^\prime\sqrt n+o(\sqrt n)$. The numerical values of
$C_1, C_2, C_3, C^\prime$ may be calculated from the law of $\bm(\cdot)$.

Going back to the observation at the end of Example~\hydrod.1, we draw the reader's
attention to the fact that the information on the law of $\bm(\cdot)$ 
that is absent did not affect the  derivation of the law of $\eta_n$ in the
present example. The reason is in that the height of the counting profile
 $g_n(\cdot)$ (i.e.
$g_n(0)$)  is indifferent to the distribution of the corresponding $\eta_n$.

The picture presented above is not the exact law of $\eta_n$ 
because $\sqrt n \bm(n^{-1}\cdot)$
might diverge from $g_n(\cdot)$  by at most $o(n^{1/2})$ along the $y$-axis 
and by at most $o(n)$ along the $x$-axis (in fact, they  must diverge since
the former does not attain its values in $\cR$, and the latter does), and
thus, (\hydrod.\true) is an approximate relation. In fact, our argument 
does not control
the particle amount with the precision  higher than $o(n^{1/2})$, that is,
the real amount of particles in a region $I$ might differ from
that of the picture presented, by at most $|I|\times o(n^{1/2})$. Also, we
were unable to give the particle distributions within each pattern, but
instead, gave its characteristics in terms of its mean values. 
These are certainly, disadvantages of the presented approach when compared to a
direct calculation of the law of $\eta_n$. An advantage is the simplicity: the
presented approach requires only to know the law of $W^{\min}(\cdot)$ while a
direct calculation requires a lot of work and produces usually cumbersome
formulae.  We note here that if one wishes to find $\eta_n$ in the exact form
then this could be done in the following way: translate CA 184 to BA 
(via $T$ from 
Assertion~\introd.\introdassone), use then the  
tools developed in [\belfer] to find the law of BA at time $n$, and, finally,
translate this law back into terms of CA 184.
\vskip4mm


\noindent{\bf \rate.\ RATE OF CONVERGENCE TO THE EQUILIBRIUM}
\vskip2mm

In this section, we present an argument that estimates both the rate of
convergence of CA 184 to its invariant state and the rate of the decay of
particle density in BA (the former is the content of
Theorem~\rate.\ratetheorem, and the latter is an auxiliary result established
in the course of the proof of this theorem). The way this argument works
is demonstrated here on a
particular example which is the CA 184 starting from the Bernoulli 1/2 product
measure. Remark~\rate.\rateremarkone\ marks the place in this argument where
main technical problems arise, when it is
adapted to other initial distributions.
\vskip2mm

\noindent{\bf Remark~\rate.\rateremarkzero.}\ \ 
Observe that that distribution at time $n$ of the CA 184 studied here, 
has been approximately found in Example~\hydrod.\hydrodexampletwo. From this
approximation, one can conclude that the probability to see two 
fixed contiguous sites of $\Z$  at the same state in CA 184 at time $n$, 
decays at  least as $const\, \times\, n^{-1/2}$. A refinement of the arguments 
from Section~\hydrod\ would reveal that $const\, \times\, n^{-1/2}$ is in
fact, the correct asymptotics for this decay. This is exactly the result of 
Theorem~\rate.\ratetheorem\ presented below. We remark here that these two
approaches have much in common, although apparently they are different. Indeed,
if one possesses a tool to estimate (\rate.\notrandomwalk) (which is the 
cornerstone of the proof of Theorem~\rate.\ratetheorem) then the same tool may
be employed to find the distribution of the process $W^{\min}(\cdot)$ which is
the basic ingredient for the approximation argument exhibited in
Section~\hydrod. 
\vskip2mm

\noindent {\bf Theorem \rate.\ratetheorem.}\ \ 
{\sl Consider CA 184 $\{\eta_n, n\in \naturals\}$ in which $\eta_0$ is
distributed by  the Bernoulli 1/2 product measure, {\it i.e.}\ $\{\eta_0(i),
i\in \Z\}$ are i.i.d. random variables with $\pr\bl \eta_0(i)=1\br=\pr\bl \eta_0(i)=0\br=1/2$. Let
$\mu_n$ denote the distribution of $\eta_n$. Then $\mu_n$ converges as 
$n\rightarrow\infty$, to  
$1/2(\delta_o+\delta_e)$ at the rate $n^{-1/2}$. In exact terms this means
that for each $k\in \naturals$ there exist $0<c^\prime_k< c_k<\infty$
such that
$$
c^\prime_k n^{-1/2}\leq \sup_{U\in \{0,1\}^{\{0,1, \ldots, k\}}}\bigl|
\mu_n(U)-1/2(\delta_o+\delta_e)(U)\bigr|\leq c_k n^{-1/2}
\eqno(\rate.\intheoremrate)
$$
}
\vskip2mm

\noindent{\it Proof.}\ \ Fix $k\in \naturals$ and define $K:=\{0,1, \ldots,
k\}$. Define
$$
E:=\{\eta\in \{0,1\}^K\, :\, \exists\, i\in \{0,1,\ldots, k-1\}{\rm\ for\ 
which\ }\eta(i)=\eta(i+1)\} 
$$
Let $\alpha, \beta\in \{0,1\}^K$ be the restrictions of respectively, $o$ and
$e$ (the configuration $e$ and $o$ have been defined in
Theorem~\invar.\invartheoremtwo)  to $K\, :\, \alpha(0):=0, \beta(0):=1$,
$\alpha(i+1):=1-\alpha(i)$, $\beta(i+1):=1-\beta(i)$, $i=0,1,\ldots,k-1$. Then
$$
E\cup \{\alpha, \beta\}=\{0,1\}^K\eqno(\rate.\tutu)
$$
We now observe that since $\mu_0$ is translation invariant then $\mu_n$ must
be such as well, and this fact leads to 
$$
\mu_n(\alpha)=\mu_n(\beta)\eqno(\rate.\titi)
$$
Since $1/2(\delta_o+\delta_e)(E)=0$ and $1/2(\delta_o+\delta_e)(\alpha)=
1/2(\delta_o+\delta_e)(\beta)=1/2$ then (\rate.\tutu) and (\rate.\titi)
imply that the supremum in (\rate.\intheoremrate) is $\mu_n(E)$. In order
to estimate this measure we shall use BA $\{\zeta_n, n\in \naturals\}$ 
such that
$$
\zeta_0=T\eta_0{\rm\ where\ }\eta_0\sim{\rm\ Bernoulli\ 1/2\ product\ measure\  and\ }T{\rm\ is\ defined\ by\ }
(\introd.\relation) 
\eqno(\rate.\initial)
$$
>From the definition of $\zeta_0$ above and Assertion~\introd.\introdassone, we
have that 
$$
\mu_n(E)=\pr\bl \, \zeta_n(i)\not=0{\rm\ for\ some\ } i\in \{1,2,\ldots, 
k\}\br
$$

Let $L=L(n)$ ($R=R(n)$) denote the event that (in the considered BA) a 
particle that originated from 
$\{-n, -n+1,\ldots, -n+k\}$ (resp., $\{n, n+1, \ldots, n+k\}$) is in the region 
$K^\prime:=\{1,2,\ldots, k\}$ at time $n$. Obviously,
$$
\pr\bl L\br\leq \mu_n(E)=\pr\bl L\cap R^c\br+\pr\bl L^c\cap R\br 
+\pr\bl L\cap R\br\leq \pr\bl L\br+\pr\bl R\br=2\pr\bl L\br
$$
whereas to derive the last inequality we used that changing $\zeta_0$ to
$-\zeta_0$ does not affect the distribution of $\zeta_0$ defined in 
(\rate.\initial). As for $\pr
\bl L\br$, one easily gets the following quite rough bounds
$$
cP\leq \pr\bl L\br\leq ckP,{\rm\ for\ an\ appropriate\ }0<c<\infty
$$
where $P$ denotes the probability to have a positive particle at site $0$ at time
$0$ and for this particle not to have survived till time $n$
 (observe, this is the same
value, if $0$ is substituted by any site of $\Z$).

To estimate the value of $P$ we  we will use the following  assertions. 
Their proofs are  straightforward.
\vskip2mm

\noindent{\bf Assertion \rate.\rateassertiontwo.}\ \  
{\sl A particle in the BA will be 
annihilated till time $n$ if and only if the distance between its initial 
position and that of its annihilating companion is $\leq 2n-1$.}
\vskip2mm

\noindent{\bf Assertion \rate.\rateassertionthree.}\ \ 
{\sl Let $\zeta_0$ be a configuration from 
$\{-1,0,1\}^{\Z}$. Assume $\zeta_0(i)=1$ for some $i$. Let $f(\cdot)$ be
the integrated profile corresponding to $\zeta_0$ and such that $f(i)=0$. 
Then $\min \{ k\geq i+2\, :\, f(k)=0\}-1$ is 
the initial position of the annihilating companion in BA $\{\zeta_n, n\in
\naturals\}$ of the particle that started from $i$.}
\vskip2mm

Let $f(\cdot)$ be the integrated process related to
$\zeta_0$ whose distribution is  (\rate.\initial) conditioned to having a
positive particle at $0$. Then due to Assertions~\rate.\rateassertiontwo\ and  
and \rate.\rateassertionthree, we conclude that 
$$
P=\pr\bl f(0)=0, f(1)>0, \ldots , f(2n-1)>0\given f(1)>0\br
\eqno(\rate.\notrandomwalk)
$$

There are several ways to get the decay rate of (\rate.\notrandomwalk). The
promptest one goes along the following lines: The random variables
$\{\zeta_0(i)\}_{i\in \Z}$ are not independent, but fall under the conditions
of Chapter~4 of [\bilin], which provide that $cn^{-1/2}f(n\cdot)\rightarrow
B(\cdot)$ for some constant $c$. From this one concludes that
(\rate.\notrandomwalk)$\sim const\times n^{-1/2}$. This 
completes the proof.\hfill$\clubsuit$
\vskip2mm

\noindent{\bf Remark~\rate.\rateremarkone.}\ \ The applicability of the
presented proof to another distribution of $\zeta_0$ (and thus, also of
$\eta_0$) depends basically on ones ability
to calculate the asymptotics of (\rate.\notrandomwalk) for this distribution. 
\vskip4mm


\noindent{\bf Acknowledgments.}\ \ The presentation of this paper benefited a
lot from discussions with J. Krug and E. Speer.
During the work on this paper, the authors were supported by grants from
CNPq, CAPES, FAPESP (Brazil), NSF (USA) and DAAD (Germany).
\vskip4mm

\noindent{\bf Bibliography.}

\vskip1mm
\vskip1mm


\item{[\belfer]} V. Belitsky, P.A. Ferrari, {\it Ballistic annihilation and
deterministic surface growth.} Journal of Statistical Physics, Vol. 80, Nos.
3/4, 1995, pp. 517-543.

\item{[\moving]} V.\ Belitsky, {\it Moving local minimum of the Brownian
motion.} \ Work in progress.

\item{[\bilin]} P. Billingsley, {\it Convergence of probability measures.}
Wiley Series in Probability and Mathematical Statistics, John Wiley\& Sons,
1968

\item{[\bennaim]} E.\ Ben-Naim, S.\ Redner, F.\ Leyvraz, {\it Decay kinetics of
ballistic annihilation.}\  Phys. Rev. Lett. {\bf 70}, 1890--1893 (1993).

\item{[\bertocci]} U. Bertocci, Surface Science Vol.15, p.296 (1969).

\item{[\bramson]} M. Bramson, C. Neuhauser, {\it Survival of one-dimensional
cellular automata under random perturbations.} Annals of Probability, Vol. 22,
No. 1, 1994, pp.244-263.

\item{[\informatic]} M. S. Capcarrere, M. Sipper, M. Tomassini, {\it
Two-state, $r=1$ cellular automaton that classifies density.} Physical Review
Letters, {\bf 77}, No. 24, 4969-4971 (1996).

\item{[\elskens]} Y.\ Elskens, H.\ L.\ Frisch, {\it Annihilation kinetics in the
one--dimensional ideal gas.}\ Physical Review A, Vol.31, No 6, June 1985.

\item{[\ermakov]} A.\ Ermakov, B. T\'oth, W.\ Werner, {\it On some annihilating
and coalescing systems.}\ \ Submitted to Journal of Statistical Physics, 1997. 

\item{[\ferravi]} P.A. Ferrari, K. Ravishankar, {\it Shocks in asymmetric
exclusion automata}. Annals of Applied Probability, Vol.2, No. 4, 1992, pp.
928-941.

\item{[\feller]} W.\ Feller, {\it Introduction to the probability theory and its
applications.} Vol. I, Wiley P, 1964.

\item{[\cyclic]} R.\ Fisch, {\it Clustering in the one--dimensional three--color
cyclic cellular automaton.} The Annals of Probability, 1992, Vol.20, No 3,
1528--1548.




\item{[\karatzas]} I.\ Karatzas, S.\ E.\ Shreve, {\it Brownian Motion and Stochastic
Calculus.} Graduate Texts in Math. Vol. 113, Springer-Verlag, 1991.

\item{[\krapivsky]} P.L. Krapivsky, F. Leyvraz, S. Redner, {\it Ballistic
Annihilation kinetics: The case of discrete velocity distribution.} Physical
Review E, Vol. 51, No. 5, 1995, pp. 3977-3987.

\item{[\krugspohn]} J.\ Krug, H.\ Spohn, {\it Universality classes for deterministic
surface growth.}\ Physical Review A, 1988, Vol. 38, 4271--4283.

\item{[\krugspohntwo]} J.\ Krug, H.\ Spohn, {\it Anomalous fluctuations in the
driven and damped sine--Gordon chain.}\ Europhys. Lett. {\bf 8} (3), 219--224
(1989). 

\item{[\liggett]} T. M. Liggett, {\it Coupling the simple exclusion process.}
Annals of Probability {\bf 4}, No. 3, 339-356 (1976).



\item{[\computersim]} T. Nagatani, {\it Creation and annihilation of traffic
jams in a stochastic asymmetric exclusion model with open boundaries: a
computer simulation.} J. Phys. A: Math. Gen. {\bf 28}, 7079-7088 (1995). 

\item{[\review]} K. Nagel, {\it Particle hopping models and traffic flow
theory.} Physical Review E, {\bf 53}, No. 5 4655-4672 (1996).

\item{[\wolfone]} S.\ Wolfram, {\it Statistical mechanics of cellular 
automata},\ Rev. Mod. Phys., 1983, Vol. 55, 601.

\item{[\wolftwo]} S.\ Wolfram, {\it Theory and Applications of Cellular
Automata.} World Scientific, Singapore, 1986.

\item{[\japonez]} H. Yaguchi, {\it Stationary measures for an exclusion 
process on one-dimensional lattices with infinitely many hopping sites.} 
Hiroshima Math. J., {\bf 16}, 449-475 (1986).
\bye